\documentclass[12pt]{article}
\usepackage{graphicx}
\usepackage[T1]{fontenc}
\usepackage{amsmath,amssymb,amsthm,microtype,tikz,booktabs,array}
\usepackage{xcolor}
\usetikzlibrary{arrows.meta,decorations.markings,positioning,calc,decorations.pathreplacing}
\usepackage[hidelinks]{hyperref}
\usepackage{float}
\usepackage[margin=31mm]{geometry}
\usepackage{hyperref,authblk}
\usetikzlibrary{shapes.geometric}
\usetikzlibrary{decorations.pathmorphing}

\newtheorem{theorem}{Theorem}
\newtheorem{lemma}[theorem]{Lemma}
\newtheorem{proposition}[theorem]{Proposition}
\newtheorem{corollary}[theorem]{Corollary}
\newtheorem{conjecture}[theorem]{Conjecture}

\newcommand{\chtt}{\chi''_{\ell}}

\title{Counterexamples to the List Total Colouring Conjecture from any regular graph}
\author[1]{W. Cames van Batenburg} 
\author[2,3]{J. Goedgebeur} 
\author[4]{G. Mazzuoccolo} 
\author[2,5]{P.~Negrini}

\affil[1]{\footnotesize D\'epartement d'Informatique, Universit\'e libre de Bruxelles, Belgium. Email: w.p.s.camesvanbatenburg@gmail.com.}
\affil[2]{\footnotesize Department of Computer Science, KU Leuven Kulak, Kortrijk, Belgium. Email: jan.goedgebeur@kuleuven.be}
\affil[3]{\footnotesize Department of Mathematics, Computer Science and Statistics, Ghent University, Ghent, Belgium.}
\affil[4]{\footnotesize Department of Physics, Computer Science and Mathematics, University of Modena and Reggio Emilia, Italy. Email: giuseppe.mazzuoccolo@unimore.it}
\affil[5]{\footnotesize Department of Mathematics and Physics, Universit\`a Cattolica del Sacro Cuore, Brescia, Italy. Email: pietro.negrini@unicatt.it}

\date{}

\begin{document}

\maketitle

\begin{abstract}
Noel recently disproved the List Total Colouring Conjecture by exhibiting a cubic graph of order $20$ with total chromatic number $4$ and list total chromatic number $5$. 
We show that such a counterexample is a particular instance of a much broader construction, which produces counterexamples from any $t$-regular graph, for every $t\geq 3$. More precisely, replacing each vertex of a $t$-regular graph $G$ with a copy of $K_{t-1,t}$ and joining corresponding terminals yields a $t$-regular graph $R_t(G)$ with total chromatic number $t+1$ and list total chromatic number $t+2$. 
Consequently, there are infinitely many regular counterexamples of every degree at least three. Moreover, they can be chosen to be bipartite.
% G: I propose "Consequently, there are infinitely many regular counterexamples of every degree at least three, even among bipartite graphs.
% W: sounds good
\end{abstract}

\section{Introduction}
Throughout this paper, all graphs are finite and simple unless otherwise stated.
A total colouring of a graph $G$ assigns colours to its vertices and edges so that adjacent vertices and edges sharing an endpoint receive different colours, and no vertex has the same colour as any of its incident edges.

%W: I read the previous variant  as forbidding that a vertex v and ALL its incident edges simultaneously have the same colour. So I've added "any of"

%A total colouring of a graph $G$ assigns colours to its vertices and
%edges so that adjacent vertices, edges sharing an endpoint, and an
%edge and either of its endpoints receive different colours. \jg{I found this a bit hard to read. What do you think about e.g.: ``A total colouring of a graph $G$ assigns colours to its vertices and
%edges so that adjacent vertices and edges sharing an endpoint receive different colours, and no vertex has the same colour as its incident edges.'' or
%``A total colouring of a graph $G$ is a colouring of its vertices and edges such that that adjacent vertices and incident edges have different colours and no vertex has the same colour as its incident edges.'' ?} {\color{red}Giuseppe: I can agree with both, but I would use adjacent edges and not incident edges} \jg{Ok! Please pick the version you like best.}

The \emph{total chromatic number} $\chi''(G)$ is the minimum number of
colours in a total colouring of $G$. Its list counterpart, the
\emph{list total chromatic number} $\chi''_{\ell}(G)$, is the least
integer $k$ such that every assignment of lists of $k$ colours to
the vertices and edges of $G$ admits a total colouring in which each
element receives a colour from its list.

Clearly, $\chi''(G)\leq\chi''_{\ell}(G)$. In the late 1990s,
Borodin, Kostochka and Woodall~\cite{BKW97}, Juvan, Mohar and
{\v{S}}krekovski~\cite{JMS98}, and Hilton and Johnson~\cite{HJ99}
independently proposed the following conjecture.

\begin{conjecture}[List Total Colouring Conjecture]
\label{conj:ltcc}
For every loopless multigraph $G$,
$$
\chi''_{\ell}(G)=\chi''(G).
$$
\end{conjecture}

Conjecture~\ref{conj:ltcc} was recently disproved by Noel~\cite{Noel}, who constructed a cubic graph $H$
on $20$ vertices satisfying
$$
\chi''(H)=4
\qquad\text{and}\qquad
\chi''_{\ell}(H)=5.
$$

His example consists of four copies of $K_{2,3}$ joined by six additional edges, one between each pair of copies, so that every vertex has degree three. In this paper, we show that Noel's construction can be generalised to produce counterexamples from any $t$-regular graph, for every integer $t\geq3$.

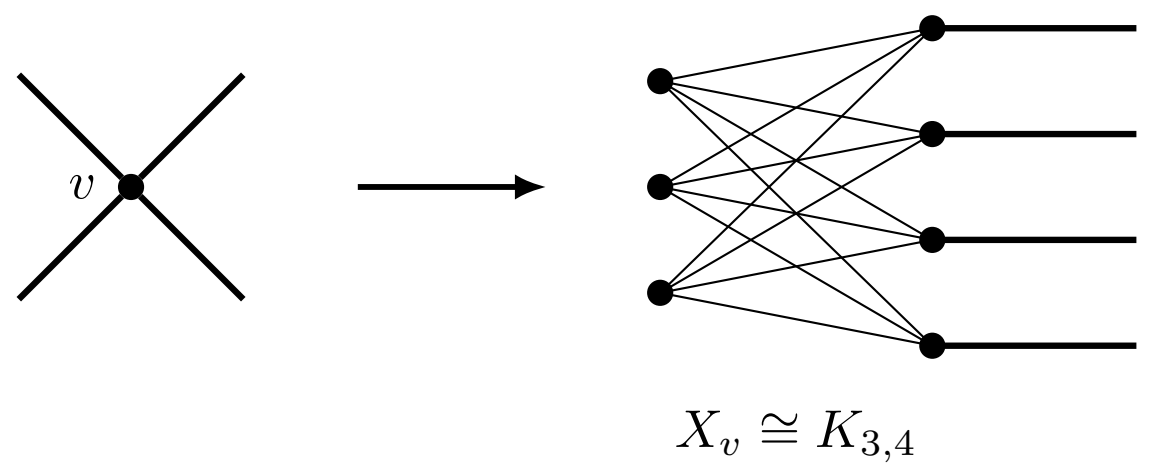
\begin{figure}[ht]
\centering
\begin{tikzpicture}[
    scale=1,
    font=\footnotesize,
    lbl/.style={font=\footnotesize, label distance=1.2pt},
    vertex/.style={circle, draw=black, fill=black, inner sep=1.6pt},
    connector/.style={very thick, black},
    intedge/.style={draw=black}
]
 \node[vertex, lbl, label=left:$v$] (v) at (0,0) {};
 \foreach \a in {45,135,225,315}{\draw[connector] (v)--(\a:1.05);}
 \draw[-{Latex[length=6pt,width=4.8pt]}, line width=1.1pt] (1.50,0) -- (2.75,0);
 \foreach \k in {0,1,2}{\coordinate (h\k) at (3.5,.7-.7*\k);}
 \foreach \j in {0,1,2,3}{\coordinate (x\j) at (5.3,1.05-.7*\j);}
 \foreach \k in {0,1,2}{\foreach \j in {0,1,2,3}{\draw[intedge] (h\k)--(x\j);}}
 \foreach \k in {0,1,2}{\node[vertex] at (h\k) {};}
 \foreach \j in {0,1,2,3}{
  \draw[connector] (x\j)--++(1.35,0);
  \node[vertex] at (x\j) {};
 }
 \node at (4.4,-1.65) {$X_v\cong K_{3,4}$};
\end{tikzpicture}
\caption{Replacing a vertex when $t=4$. In general the block is
$K_{t-1,t}$, with one connector at each of its $t$ terminals.}
\label{fig:replacement}
\end{figure}

Let $G$ be a 
%connected
% G: we should omit also simple. We declare it at the beginning of the paper.
$t$-regular graph, where $t\ge3$. Replace each
vertex $v$ by a disjoint copy $X_v$ of the complete bipartite graph $K_{t-1,t}$  (see Figure~\ref{fig:replacement}). We refer to $X_v$ as the \emph{block} associated with $v$; we call the vertices in the smallest part of $X_v$ \emph{private vertices} and those in its largest part
\emph{terminals}. Associate the terminals of $X_v$ bijectively with the
edges of $G$ incident to $v$. For each edge $uv\in E(G)$, join the
corresponding terminals of $X_u$ and $X_v$ by an edge, called a
\emph{connector}. Denote the resulting graph by $R_t(G)$.
%\textcolor{red}{W: could consider to drop the subscript $t$ in $R_t(G)$, because the size of the blow-up already follows from the regularity of $G$. On the other hand, the subscript serves as a useful recurring reminder for the reader.}.

%{\color{blue}P: Personally, I feel that especially in formulas like $\chi''(R_t(G))=t+1$, having $t$ explicit on both sides makes the dependence immediately clear to the reader}
% W: agreed. Let's keep it like this.
% G: I also agree.

The graph $R_t(G)$ is 
%connected and 
$t$-regular and has $(2t-1)|V(G)|$ vertices. Noel's example is recovered by taking $t=3$ and $G=K_4$.

The following theorem is the main result of the paper. 

\begin{theorem}\label{thm:main}
For every $t$-regular graph $G$, with $t\ge3$,
$$\chi''(R_t(G))=t+1
\qquad\text{and}\qquad
\chtt(R_t(G))= t+2.$$
\end{theorem}

Since $G$ can be arbitrarily chosen, our construction yields infinitely many counterexamples of each regular degree at least three, all with a gap of exactly one between the two parameters.
Note also that $R_t(G)$ is bipartite if and only if $G$ is, so we obtain bipartite counterexamples with an arbitrarily large total chromatic number. Moreover, every $k$-connected  choice of $G$ yields a $k$-connected counterexample, for every $k\le t$.
%In particular, by choosing $G$ bipartite, we obtain bipartite counterexamples with an arbitrarily large total chromatic number.
%Also note that every $k$-connected  choice of $G$ yields a $k$-connected counterexample, for every $k\le t$. 
%Although the construction is described for connected graphs, all our results also hold for disconnected graphs, since the construction and both colouring parameters behave componentwise.

The proof of Theorem~\ref{thm:main} is divided into three parts. In Section~\ref{sec:totalcoloring}, we establish Proposition~\ref{prop:total}, showing that $\chi''(R_t(G))=t+1$. Section~\ref{sec:lowerbound} contains Proposition~\ref{thm:lower}, which establishes the lower bound $\chtt(R_t(G)) \geq t+2$, the key step in refuting the conjecture. Lastly, Proposition~\ref{thm:upper} in Section~\ref{sec:upperbound} completes the proof of Theorem~\ref{thm:main} by showing that $\chtt(R_t(G)) \leq t+2$. We then conclude in Section~\ref{sec:final} with some final remarks.

We stress that only the results in Sections~\ref{sec:totalcoloring} and~\ref{sec:lowerbound} are necessary to establish that our constructions are counterexamples to the List Total Colouring Conjecture, while Section~\ref{sec:upperbound} shows that they cannot provide an even larger gap between $\chi''$ and~$\chtt$.

\section{A total colouring of $R_t(G)$ with $t+1$ colours}\label{sec:totalcoloring}

\begin{proposition}\label{prop:total}
For every $t$-regular graph $G$, with $t\geq 3$,
$$\chi''(R_t(G))=t+1.$$
\end{proposition}
\begin{proof}
By the $t$-regularity of $R_t(G)$, any total colouring requires at least $t+1$ colours. Hence, to obtain $\chi''(R_t(G))= t+1$, we only have to prove the existence of a total colouring with $t+1$ colours.

Let $J_t(G)$ be the graph obtained from $R_t(G)$ by deleting all private vertices and turning the terminals of each block into a clique. 
Since $J_t(G)$ is $t$-regular and no connected component is a complete graph (nor an odd cycle, since $t \geq 3$), Brooks' theorem~\cite{RLB} guarantees a proper vertex $t$-colouring $\alpha$, using colours $0,1,2,\ldots, t-1$. Since the vertices corresponding to the terminals in the same block of $R_t(G)$ form a clique in $J_t(G)$, they receive all $t$ distinct colours. 

Now keep these colours on the corresponding vertices in $R_t(G)$. Consider a block and let $u_1,\ldots,u_{t-1}$ be its private vertices. Assign the colour $(x+j)\bmod t$ to the internal edge between the private vertex $u_j$ and the terminal with colour $x$.
Hence, at a terminal of colour $x$, the incident internal edges use precisely the $t-1$ colours other than $x$, while at a private vertex they use all $t$ colours. Apply the same colouring rule to every block of $R_t(G).$

Finally, assign a new colour $\infty$ to every private vertex and every connector. Since the connectors form a matching and $\alpha$ assigns different colours to their endpoints, this gives a proper total colouring of $R_t(G)$ with $t+1$ colours.
\end{proof}

\section{A lower bound of $t+2$ for $\chtt(R_t(G))$}\label{sec:lowerbound}

To prove that $\chtt(R_t(G))\geq t+2$, it suffices to exhibit a list assignment of size $t+1$ to the vertices and edges of $R_t(G)$ that admits no proper total colouring.

We first establish a simple property of total colourings of $K_{t-1,t}$ with $t+1$ colours.

\begin{lemma}\label{lem:block}
Let $A$ and $B$ be the two parts of $K_{t-1,t}$, with $|A|=t-1$ and $|B|=t$. Let $L$ be a set of $t+1$ colours. In every proper total colouring of $K_{t-1,t}$ using colours from $L$:
\begin{itemize}
\item the vertices of $B$ receive pairwise distinct colours;
\item all vertices of $A$ receive the same colour $\ell\in L$;
\item for each vertex $b\in B$, the colours assigned to $b$ and its incident edges are precisely those of $L\setminus\{\ell\}$.
\end{itemize}

\end{lemma}
\begin{proof}
  Suppose, by contradiction, that two vertices of $B$ have the same colour $c$. Hence, no vertex in $A$ can have colour $c$. A vertex of $A$ and its $t$ incident edges use all $t+1$ colours of $L$, so each vertex of $A$ has an incident edge coloured $c$. These $t-1$ edges must meet distinct vertices in $B$ and cannot meet the two vertices with colour $c$. This requires $t-1$ distinct vertices among only $t-2$ available, a contradiction.

Since every vertex of $A$ is adjacent to all vertices of $B$, all vertices of $A$ must receive the remaining colour, say $\ell$. Moreover, neither a vertex of $B$ nor an edge can have colour $\ell$. Hence, the colours used by a vertex of $B$ and its $t-1$ incident edges are the $t$ colours of $L\setminus\{\ell\}$.
\end{proof}

When applying Lemma~\ref{lem:block} to a block of $R_t(G)$, we say that the block \textit{has list $L$} if all of its vertices and internal edges are assigned the same list $L$. In such a colouring, we call the common colour $\ell$ of the private vertices the \textit{private colour} of the block.

\begin{corollary}\label{cor:privatecolor}
    If a block of $R_t(G)$ has list $L$ of size $t+1$ and private colour $\ell$, then no incident connector can receive a colour from $L\setminus{\{\ell\}}$. Moreover, every incident connector with list $L$ receives colour $\ell$.
\end{corollary}

\begin{proposition}\label{thm:lower}
For every $t$-regular graph $G$, with $t\ge3$,
$$\chi''_{\ell}(R_t(G))\geq t+2.$$
\end{proposition}
\begin{proof}

Set $a=\lfloor(t+1)/2\rfloor$ and $b \in \{0,1\}$ such that $b \not\equiv t \pmod 2$. Choose pairwise disjoint colour sets $A_1,A_2,A_3,E$, such that $|A_i|=a$ for each $i$ and $|E|=b$. We define three lists of size $t+1$:
\[ L_1=A_2\cup A_3\cup E, \quad L_2=A_1\cup A_3\cup E \quad \text{and} \quad L_3=A_1\cup A_2\cup E.
\]

Choose a set $L_4\subseteq A_1\cup A_2\cup A_3$ of size $t+1$, which is possible since $3a\geq t+1$.

Since $G$ is $t$-regular with $t\geq 3$, it contains a subdivision $T$ of $K_4$ (see~\cite{Dirac}). Let $v_1,v_2,v_3,v_4$ be its branch vertices, and, for $i<j$, let $p_{ij}$ denote the path in $T$ joining $v_i$ and $v_j$ with no other branch vertices. Write $X_i=X_{v_i}$, and let $P_{ij}$ denote the chain of blocks consisting of the connectors corresponding to the edges of $p_{ij}$ and the blocks corresponding to its internal vertices and thus excluding $X_i$ and $X_j$.

We say that a chain of blocks $P_{ij}$ has list $L$ if all of its connectors and all vertices and internal edges of its blocks are assigned list $L$. We assign lists to the blocks $X_1,X_2,X_3,X_4$ and the chains $P_{ij}$ as follows (see Figure~\ref{fig:list-assignment}):
\begin{itemize}
    \item assign the list $L_i$ to the block $X_i$, for $i=1,2,3,4$;
    \item assign the list $L_4$ to each of the chains $P_{14},P_{24},P_{34}$, and lists $L_1,L_2,L_3$ to $P_{23},P_{13},P_{12}$, respectively.
\end{itemize}

\begin{figure}[ht]
\centering

\begin{tikzpicture}[
    block/.style={
        circle,
        draw,
        line width=0.6pt,
        minimum size=8mm,
        inner sep=1pt,
        align=center
    },
    chain/.style={
        draw,
        line width=0.6pt,
        decorate,
        decoration={
            snake,
            amplitude=0.5mm,
            segment length=0.8mm
        }
    },
    lab/.style={
        font=\small,
        fill=white,
        inner sep=1.5pt
    }
]

% Equilateral triangles
\pgfmathsetmacro{\h}{2.2*sqrt(3)}
\pgfmathsetmacro{\g}{\h/3}

% ==========================================
% LEFT: Distinguished blocks and chains
% ==========================================

\begin{scope}[xshift=-2.8cm]

\node[block] (X1) at (0,\h) {$X_1$};
\node[block] (X2) at (-2.2,0) {$X_2$};
\node[block] (X3) at (2.2,0) {$X_3$};
\node[block] (X4) at (0,\g) {$X_4$};

% Outer chains

\draw[chain] (X1) --
    node[midway,lab] {$P_{12}$} (X2);

\draw[chain] (X1) --
    node[midway,lab] {$P_{13}$} (X3);

\draw[chain] (X2) --
    node[midway,lab] {$P_{23}$} (X3);

% Inner chains

\draw[chain] (X1) --
    node[midway,lab] {$P_{14}$} (X4);

\draw[chain] (X2) --
    node[midway,lab] {$P_{24}$} (X4);

\draw[chain] (X3) --
    node[midway,lab] {$P_{34}$} (X4);

\end{scope}

% ==========================================
% RIGHT: List assignment
% ==========================================

\begin{scope}[xshift=2.8cm]

\node[block] (Y1) at (0,\h) {$\boldsymbol{L_1}$};
\node[block] (Y2) at (-2.2,0) {$\boldsymbol{L_2}$};
\node[block] (Y3) at (2.2,0) {$\boldsymbol{L_3}$};
\node[block] (Y4) at (0,\g) {$\boldsymbol{L_4}$};

% Outer chains

\draw[chain] (Y1) --
    node[midway,lab] {$\boldsymbol{L_3}$} (Y2);

\draw[chain] (Y1) --
    node[midway,lab] {$\boldsymbol{L_2}$} (Y3);

\draw[chain] (Y2) --
    node[midway,lab] {$\boldsymbol{L_1}$} (Y3);

% Inner chains

\draw[chain] (Y1) --
    node[midway,lab] {$\boldsymbol{L_4}$} (Y4);

\draw[chain] (Y2) --
    node[midway,lab] {$\boldsymbol{L_4}$} (Y4);

\draw[chain] (Y3) --
    node[midway,lab] {$\boldsymbol{L_4}$} (Y4);

\end{scope}

\end{tikzpicture}

\caption{Schematic representation of the distinguished
blocks $X_i$ and chains $P_{ij}$ (left), together with
their assigned lists (right), as used in
Proposition~\ref{thm:lower}.}

\label{fig:list-assignment}

\end{figure}
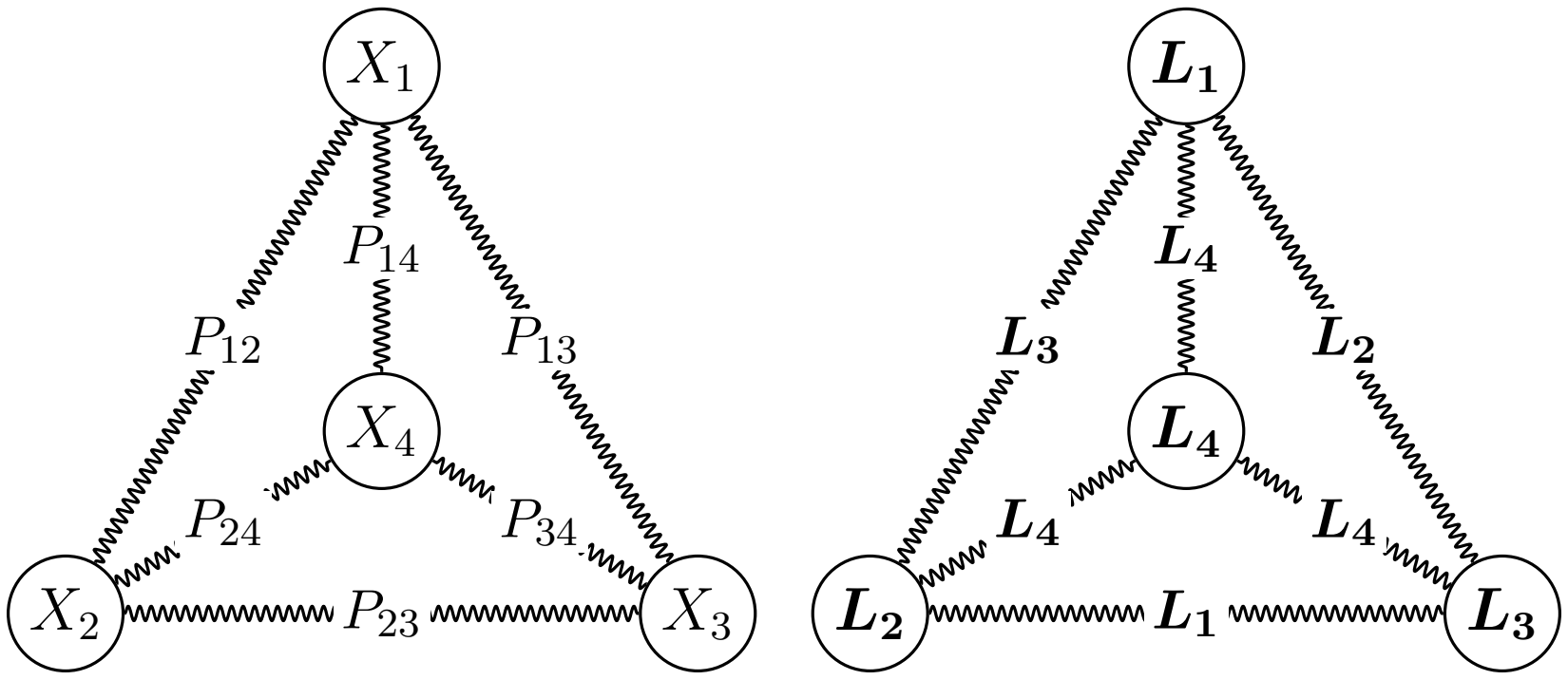

Suppose, for a contradiction, that there exists a total colouring respecting these lists. All connectors of the chain $P_{ij}$ must receive the same colour: at each internal block, the two connectors of the chain share the block's list, so both receive its private colour by Corollary~\ref{cor:privatecolor}. If $p_{ij}$ has no internal vertices, $P_{ij}$ consists of a single connector.

Let $c\in L_4$ be the private colour of $X_4$.
Since the chains $P_{14},P_{24},P_{34}$ all have list $L_4$, their connectors incident to $X_4$, by Corollary~\ref{cor:privatecolor}, must receive colour $c$. By the previous observation, every connector in these three chains must have colour $c$. 
Since $L_4\subseteq A_1\cup A_2\cup A_3$, there is a unique $i\in\{1,2,3\}$ such that $c\in A_i$. Without loss of generality, say $i=1$.

Since both $L_2$ and $L_3$ contain $c$, the connectors of colour $c$
incident to $X_2$ and $X_3$ force both blocks to have private colour $c$.
On the other hand, all connectors of $P_{23}$ must receive a common colour $d \in L_1=A_2\cup A_3\cup E.$
Since $c \notin L_1$, $d\ne c$.
But $L_1\subseteq L_2\cup L_3$, so $d$ belongs to at least one
of the two block lists of $X_2$ and $X_3$.
At that endpoint, Corollary~\ref{cor:privatecolor} forces the connector's colour to equal the block's private colour $c$, contradicting $d\ne c$.
Thus this $(t+1)$-list-assignment admits no total colouring.
Extending the prescribed lists arbitrarily to the remaining vertices and edges of $R_t(G)$ gives an uncolourable $(t+1)$-list assignment, completing the proof.
\end{proof}

\section{A matching upper bound for $\chtt(R_t(G))$}\label{sec:upperbound}

For a vertex list assignment
$L=\{L(v)\}_{v\in V(G)}$, we say that $G$ is \textit{$L$-colourable} if it admits a proper vertex colouring $c$ such that $c(v)\in L(v)$ for every $v\in V(G)$. Similarly, for a total list assignment $L=\{L(x)\}_{x\in V(G)\cup E(G)}$, we say that $G$ is \textit{totally $L$-colourable} if it admits a proper total colouring $c$ such that $c(x)\in L(x)$ for every $x\in V(G)\cup E(G)$.

To complete the proof of Theorem~\ref{thm:main}, it remains to show that $\chi''_{\ell}(R_t(G))\le t+2$.

\begin{proposition}\label{thm:upper}
For every $t$-regular graph $G$, with $t\ge3$,
$$\chi''_{\ell}(R_t(G))\le t+2.$$
\end{proposition}

We first reduce the proposition to a local assertion. 

Let
$A=\{a_1,\ldots,a_{t-1}\}$ and $B=\{b_1,\ldots,b_t\}$ be the two parts of
$K_{t-1,t}$, and denote its edges by $e_{i,j}=a_ib_j$.
We will show that if every total list assignment satisfying
$$ |L(a_i)|\ge t+2,\qquad |L(b_j)|\ge t,\qquad |L(e_{i,j})|\ge t+1
$$
admits a total colouring of $K_{t-1,t}$, then $R_t(G)$ admits a total colouring from lists of at least $t+2$ colours. 

Start with lists of at least $t+2$ colours on the vertices and edges of $R_t(G)$ and colour the connectors arbitrarily from their lists (which causes no conflicts, as the connectors form a matching).
Process the blocks in any order. At each terminal, delete its connector's
colour from its list if present, and if the terminal at the other end of that connector has already been coloured, delete its colour as well where applicable. At an edge inside the block, remove the adjacent connector's colour contained in its list. The lists of private vertices do not lose colours. The remaining list sizes are at least $t+1$ for the edge $e_{i,j}$, at least $t$ for a terminal $b_j$, and $t+2$ for a private vertex $a_i$. Apply the local assertion to colour the current block from these lists, and continue with the next block. These deletions also resolve every conflict that involves different blocks or connectors.

Thus, it only remains to prove the following lemma.

\begin{lemma} \label{lemma:winning-list}
    Let $t \geq 3$. For any total list assignment $L$ of $K_{t-1,t}$, with $|L(a_i)| \geq t+2$, $|L(b_j)| \geq t$ and $|L(e_{i,j})| \geq t+1$, the graph $K_{t-1,t}$ is totally $ L$-colourable.
\end{lemma}

We now recall the tools needed to prove Lemma~\ref{lemma:winning-list}.\\

Let $G=(V,E)$ be a graph and $\rho : V \to \mathbb{N}$. An orientation $D$ of $G$ is called a \textit{$\rho$-orientation} if $d_D^+(v) = \rho(v)$ for every vertex $v \in V$. Fix a reference orientation $D_0$ of $G$. We say that an orientation $D$ is \textit{even} if it has an even number of edges oriented in reverse with respect to $D_0$, otherwise we will call it \textit{odd}.

\begin{theorem}[Alon and Tarsi~\cite{AT}] \label{theorem:alon-tarsi}
    Let $G$ be a graph and $\rho : V \to \mathbb{N}$. If the number of even $\rho$-orientations is not equal to the number of odd $\rho$-orientations, then for every vertex list assignment $L=\{L(v)\}_{v \in V}$ satisfying $|L(v)| \geq \rho(v)+1$, $G$ is $ L$-colourable.
\end{theorem}

A \textit{clique decomposition} of a graph $G$ is a partition of the edge set of $G$ into edge disjoint cliques $G_1, \dots, G_\ell$, denoted by $G_1 \oplus \dots \oplus G_\ell$. A \textit{clique transitive orientation} of $G$ is an orientation whose restriction to each clique of the specified decomposition is a transitive tournament.

The following lemma shows that to verify the condition in the Alon-Tarsi theorem, it is sufficient to consider clique transitive orientations.

\begin{lemma}[H{\"a}ggkvist and Janssen~\cite{HJ}] \label{lemma:clique-dec}
    Let $G = G_1 \oplus \dots \oplus G_\ell$ be a clique decomposition of $G$ and $\rho: V \to \mathbb{N}$. The number of even $\rho$-orientations that are not clique transitive is equal to the number of odd $\rho$-orientations that are not clique transitive.
\end{lemma}

In particular, if there is exactly one clique transitive $\rho$-orientation, then the numbers of even and odd $\rho$-orientations differ by one and the condition of the Alon-Tarsi theorem holds.

As a consequence we obtain the following corollary.

\begin{corollary} \label{corollary:unique-decomposition}
    Let $G = G_1 \oplus \dots \oplus G_\ell$ be a clique decomposition of $G$, let $\rho: V(G) \to \mathbb{N}$, and let $L=\{L(v)\}_{v \in V}$ be a vertex list assignment  satisfying  $|L(v)| \geq \rho(v)+1$. If there exists a unique clique transitive $\rho$-orientation $D$, then $G$ is $ L$-colourable.
\end{corollary}

Let us now consider the total graph $T(K_{t-1,t})$. Recall that the vertex set $V(T(K_{t-1,t})) = A \cup B \cup \{ e_{i,j}:i \in \{1, \dots, t-1\}, j \in \{1, \dots, t\} \}$, where we use the same symbol $e_{i,j}$ for the edge of $K_{t-1,t}$ and its corresponding vertex in $T(K_{t-1,t})$.  We also identify a total list assignment of $K_{t-1,t}$ with the corresponding vertex list assignment of $T(K_{t-1,t})$.

Following the approach of H{\"a}ggkvist and Janssen~\cite{HJ}, we consider an augmented graph $T^*(K_{t-1,t})$ obtained from $T(K_{t-1,t})$ by adding $t-1$ new vertices $u_1, \dots, u_{t-1}$, and joining each $u_i$ to $a_i$ and to all $e_{i,1}, \dots, e_{i,t}$. In $T^*(K_{t-1,t})$ we define the following cliques:
\begin{itemize}
    \item the row clique $R_i$, induced by $u_i$, $a_i$, $e_{i,1}, \dots, e_{i,t}$, of size $t+2$;
    \item the column clique $C_j$, induced by $b_j$, $e_{1,j}, \dots, e_{t-1,j}$, of size $t$;
    \item the $2$-clique $Q_{i,j}$, on the vertices $a_i$ and $b_j$.
\end{itemize}
These cliques give a clique decomposition of $T^*(K_{t-1,t}) = R_1 \oplus \dots \oplus R_{t-1} \oplus C_1 \oplus \dots \oplus C_{t} \oplus Q_{1,1} \oplus \dots \oplus Q_{t-1,t}$. To apply Corollary~\ref{corollary:unique-decomposition}, it suffices to establish uniqueness
with respect to this fixed clique decomposition.

In what follows, we use the notation $u \to v$ to indicate that the edge $uv$ is oriented from $u$ to $v$.

\begin{proposition} \label{proposition:augmented}
    For $i \in \{1, \dots, t-1\}$ and $j \in \{1, \dots, t\}$, let $\rho: V(T^*(K_{t-1,t})) \to \mathbb{N}$ be such that $\rho(u_i) =t$, $\rho(b_j) = t-1$, $\rho(a_i) = t+1$ and  
    $$
    \rho(e_{i,j}) =
    \begin{cases} 
        t & \text{ if } \quad 2 \leq j \leq t+1-i\\
        t-1 & \text{ if } \quad j=1 \quad \text{or} \quad j > t+1-i.   
    \end{cases}
    $$
    Then, $T^*(K_{t-1,t})$ has a unique clique transitive $\rho$-orientation with respect to the clique decomposition into the row cliques $R_i$, the column cliques $C_j$, and the $2$-cliques $Q_{i,j}$ defined above.
\end{proposition}

\begin{proof}
We first show that there is at most one such orientation.

    \begin{itemize}
        \item Consider a $(t+2)$-clique $R_i$. Since $R_i$ is transitively oriented and has order $t+2$, it contains a vertex of local out-degree $t+1$. This vertex must be $a_i$, since every other vertex of $R_i$ has total out-degree at most $t$.
        Moreover, $u_i$ has all its neighbours in $R_i$; recalling that  $a_i  \to u_i$ and $\rho(u_i)=t$, we have that  $u_i \to e_{i,j} $ holds for all $j=1,\ldots,t$.
        \item Consider a $2$-clique $Q_{i,j}$. Since $\rho(a_i)=t+1$, and $a_i$ already has out-degree $t+1$ in $R_i$, all remaining edges incident to $a_i$ must be directed to $a_i$, in particular $b_j \to a_i$ in each $Q_{i,j}$.
        \item Consider a $t$-clique $C_j$. Since $\rho(b_j)=t-1$ and $b_j \to a_i$, all of the remaining edges incident to $b_j$ must be directed to $b_j$, in particular $e_{i,j} \to b_j$ for all $i$.
    \end{itemize}

Figure~\ref{fig:k34-R2-C2-b3} illustrates some of these forced orientations in the case $t=4.$

\begin{figure}[H]
\centering
\begin{tikzpicture}[
scale=0.8,
    vtx/.style={circle, draw=black, fill=black, inner sep=1.6pt},
    lbl/.style={font=\footnotesize, label distance=1.2pt},
    arr/.style={
    thick,
    draw=gray!60,
    postaction={decorate},
    decoration={
        markings,
        mark=at position 0.5 with {
            \arrow[black]{Latex[length=4.5pt,width=3.5pt]}
        }
    }
}
]    
\def\d{1.8}

\foreach \i in {1,2,3}{%
    \pgfmathsetmacro{\y}{(3-\i)*\d}
    \node[vtx, lbl, label=above left:$u_{\i}$] (u\i) at (-2*\d, \y) {};
    \node[vtx, lbl, label=above left:$a_{\i}$] (a\i) at (-\d, \y) {};
}
\foreach \j in {1,2,3,4}{%
    \pgfmathsetmacro{\x}{(\j-1)*\d}
    \node[vtx, lbl, label=above:$b_{\j}$] (b\j) at (\x, 3*\d) {};
}

\foreach \i in {1,2,3}{%
    \pgfmathsetmacro{\y}{(3-\i)*\d}
    \foreach \j in {1,2,3,4}{%
        \pgfmathsetmacro{\x}{(\j-1)*\d}
        \node[vtx, lbl, label=right:$e_{\i,\j}$] (e\i\j) at (\x, \y) {};
    }
}
\draw[arr] (a2) -- (u2);
\draw[arr] (a2) -- (e21);
\draw[arr] (a2) to[bend left=16] (e22);
\draw[arr] (a2) to[bend left=20] (e23);
\draw[arr] (a2) to[bend left=24] (e24);

\draw[arr] (u2) to[bend right=18] (e21);
\draw[arr] (u2) to[bend right=21] (e22);
\draw[arr] (u2) to[bend right=24] (e23);
\draw[arr] (u2) to[bend right=27] (e24);

\draw[arr] (e12) -- (b2);
\draw[arr] (e22) to[bend left=22] (b2);
\draw[arr] (e32) to[bend left=26] (b2);

\draw[arr] (b3) .. controls (1.5, 5.10) and (-0.5, 4.75) .. (a1);
\draw[arr] (b3) .. controls (1.35, 4.95) and (-0.8, 4.50) .. (a2);
\draw[arr] (b3) .. controls (1.20, 4.85) and (-1.1, 4.30) .. (a3);

\end{tikzpicture}
\caption{The $22$ vertices of $T^*(K_{3,4})$, showing the orientations of the edges of $R_2$ incident to $a_2$ and $u_2$, the edges of $C_2$ incident to $b_2$, and the edges $b_3 \to a_i$ ($i=1,2,3$).}
\label{fig:k34-R2-C2-b3}
\end{figure}
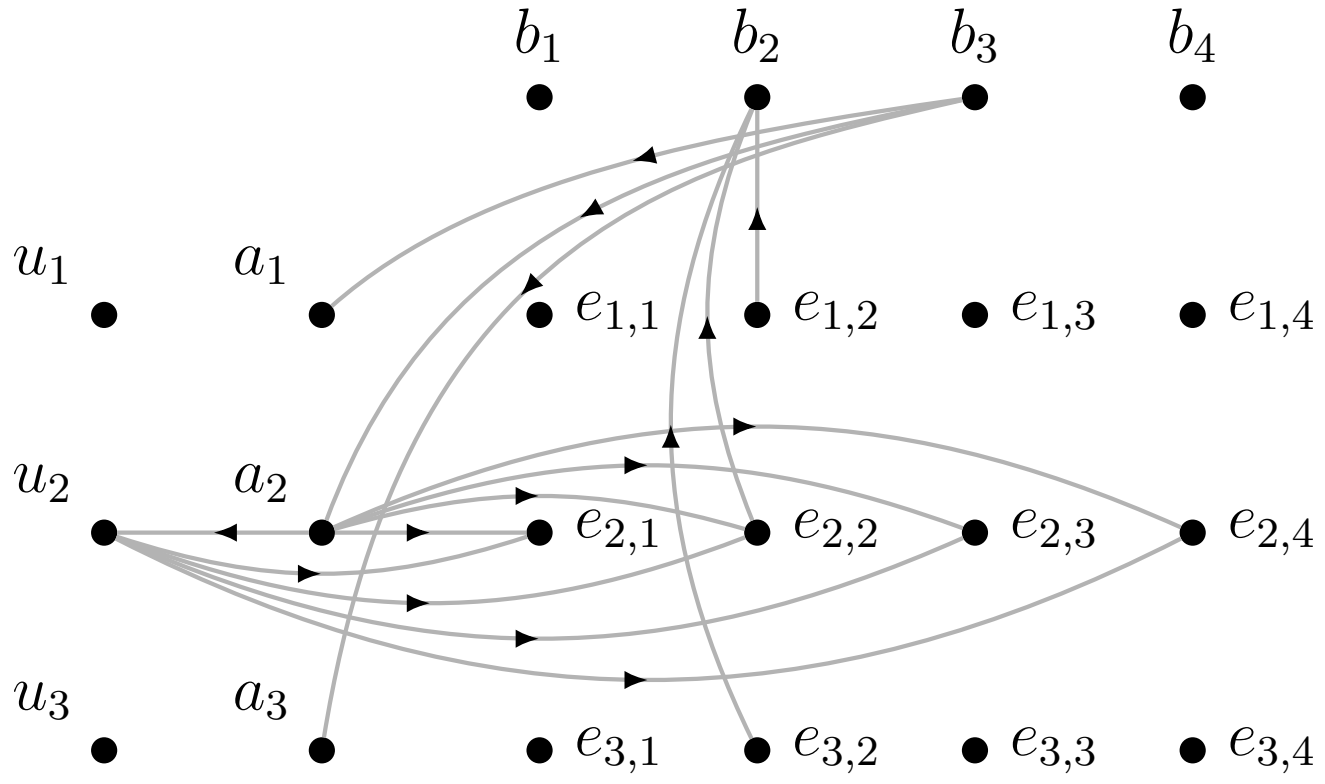

The only edges whose orientations remain to be determined
are those of the line graph $L(K_{t-1,t})$.
Each vertex $e_{i,j}$ has exactly one outgoing edge outside
the line graph, namely $e_{i,j}\to b_j$.
Hence, its required out-degree within $L(K_{t-1,t})$ is

    $$
        \rho(e_{i,j})-1=
        \begin{cases}
            t-1 & \text{ if } \quad 2 \leq j \leq t+1-i\\
            t-2 & \text{ if } \quad j=1 \quad \text{or} \quad j > t+1-i.
        \end{cases}
    $$
    
By Janssen's Lemma 2.4 of~\cite{JJ}, applied with parameters $n=t$ and $r=t-1$, there exists a unique clique transitive orientation of $L(K_{t-1,t})$ that has the prescribed  out-degrees. 

Conversely, Janssen's orientation can be extended
by orienting all remaining edges as described above. Each row clique is transitive with $a_i$ first and $u_i$ second, and each column clique is transitive with $b_j$ last. This gives a clique transitive $\rho$-orientation. This completes the proof.
\end{proof}

\begin{proof}[Proof of Lemma~\ref{lemma:winning-list}]
Let $ L$ be a total list assignment of $K_{t-1,t}$
satisfying the hypotheses of the lemma.
We regard $ L$ as a vertex list assignment of the
total graph $T(K_{t-1,t})$, identifying each edge of
$K_{t-1,t}$ with its corresponding vertex in $T(K_{t-1,t})$.

Extend $ L$ to a vertex list assignment $ L'$
of $T^*(K_{t-1,t})$ by setting $L'(x)=L(x)$ for every
$x\in V(T(K_{t-1,t}))$ and assigning an arbitrary list
of $t+1$ colours to each added vertex $u_i$.

By Proposition~\ref{proposition:augmented}, the graph $T^*(K_{t-1,t})$ has a unique
clique transitive $\rho$-orientation with respect to the
specified clique decomposition. Corollary~\ref{corollary:unique-decomposition} therefore
implies that $T^*(K_{t-1,t})$ admits an $ L'$-colouring.

Since $ L'$ agrees with $ L$ on $V(T(K_{t-1,t}))$, restricting this colouring to $T(K_{t-1,t})$ gives an $ L$-colouring of the total graph. Hence $K_{t-1,t}$ is totally $ L$-colourable. 
\end{proof}

Together with the preceding reduction, this completes the proof of Proposition~\ref{thm:upper}.

\section{Final remarks}\label{sec:final}
%\begin{itemize}
%\item
Although our construction produces regular graphs, the proof of Proposition~\ref{thm:lower} also yields non-regular counterexamples. Indeed, let $F_T$ be the subgraph of $R_t(G)$ consisting of the blocks corresponding to the vertices of the subdivision $T$, together with the connectors corresponding to its edges. The proof gives $\chi''_\ell(F_T)\geq t+2$. On the other hand, $F_T$ has maximum degree $t$, and the total colouring from Proposition~\ref{prop:total} restricts to $F_T$. Thus $\chi''(F_T)=t+1$, and every non-regular subgraph of $R_t(G)$ containing $F_T$ is also a counterexample.
%\item

More generally, the list lower bound in  Proposition~\ref{thm:lower} remains valid if only the vertices of $T$, or any larger set, are replaced in $G$. However, the resulting $t$-regular graph need not admit a total $(t+1)$-colouring. Indeed, if the replaced vertices form a connected subgraph and some unreplaced vertex has at least two neighbours in it, Lemma~\ref{lem:block} forces the corresponding connectors to receive the same colour, a contradiction.
%\item

Woodall~\cite{Woodall2006} and Hetherington–Woodall~\cite{HetheringtonWoodall2006} showed that the List Total Colouring Conjecture does hold for all $K_4$-minor-free graphs.
Thus, every counterexample must contain a subdivision of $K_4$.
Interestingly, our construction for $R_t(G)$ exploits exactly such a structure, with the lower bound argument in Proposition~\ref{thm:lower} relying on a subdivision of $K_4$ in the underlying graph $G$. In fact, our main result extends to any loopless $t$-regular multigraph $G$ containing a subdivision of $K_4$. 
%\item

Finally, as already remarked by Noel in~\cite{Noel}, it remains an interesting problem to determine whether the gap between the list total chromatic number and the total chromatic number can be at least two. We also recall that a gap of at least three would disprove the List Edge-Colouring Conjecture, since $\chi''_\ell(G)\leq\chi'_\ell(G)+2$ and $\chi'(G)\leq\chi''(G)$.
The latter conjecture has been attributed to several researchers; for an account of its early history, see~\cite[Problem~12.20]{JT}.
%\end{itemize}

\section*{Acknowledgements}
We thank Jonathan Noel for discussions.
Wouter Cames van Batenburg is supported by the Belgian National Fund for Scientific Research (FNRS). Jan Goedgebeur is supported by Internal Funds of KU Leuven and a grant of the Research Foundation Flanders (FWO) with grant number G0AGX24N. Giuseppe Mazzuoccolo and Pietro Negrini are partially supported by a grant of the group GNSAGA of INdAM.

\section*{Declaration of use of AI}
A closed-weight AI model was used for general language editing, initial versions of the figures, and to identify the key step of Section~\ref{sec:lowerbound} to use a subdivision of $K_4$.

\bibliographystyle{plain}
\bibliography{references}

\end{document}